\documentclass{autart}

\usepackage{graphics} 
\usepackage{epsfig} 
\usepackage{mathptmx} 
\usepackage{times} 
\usepackage{amsmath} 
\usepackage{amssymb} 
\usepackage{amsfonts}
\usepackage{psfrag}
\usepackage{graphicx}
\usepackage[comma]{natbib}
\usepackage{color}

\newtheorem{theorem}{Theorem}[section]

\def\bi{\begin{itemize}}
\def\ei{\end{itemize}}
\def\bn{\begin{enumerate}}
\def\en{\end{enumerate}}
\def\bq{\begin{eqnarray}}
\def\eq{\end{eqnarray}}
\def\be{\begin{equation}}
\def\ee{\end{equation}}
\def\bea{\begin{eqnarray}}
\def\eea{\end{eqnarray}}
\def\beann{\begin{eqnarray*}}
\def\eeann{\end{eqnarray*}}
\def\bsea{\begin{subeqnarray}}
\def\esea{\end{subeqnarray}}
\def\bmat{\left[ \begin{array}}
\def\emat{\end{array} \right]}

\def\ns{\hspace{-1mm}}
\def\proof{\noindent{\bf{\em Proof:}\ \ }}
\def\QED{\mbox{\rule[0pt]{1.5ex}{1.5ex}}}
\def\endproof{\hspace*{\fill}~\QED\par\endtrivlist\unskip}
\newcommand{\real}{{\mathbb{R}}}

\def\spacingset#1{\def\baselinestretch{#1}\small\normalsize}
\spacingset{0.97}

\def\gL{{\mathcal L}}
\def\gM{{\mathcal M}}

\def\gP{{\mathcal P}}

\newfont{\BB}{msbm10}
\newfont{\bb}{msbm8}

\newcommand{\bmx}{\begin{matrix}}
\newcommand{\emx}{\end{matrix}}
\newcommand{\ba}{\begin{array}}
\newcommand{\ea}{\end{array}}
\newcommand{\defi}{\stackrel{\text{\tiny def}}{=}}
\def\undp{{\underline \pi}}
\def\overp{{\overline \pi}}

\newcommand{\rank}{\operatorname{rank}}

\newcommand{\diag}{\operatorname{diag}}
\newcommand{\blkdiag}{\operatorname{blkdiag}}
\def\nn{\nonumber}
\def\bq{\begin{eqnarray}}
\def\eq{\end{eqnarray}}
\def\tra{{\scalebox{.6}{\thinspace\mbox{T}}}}
\def\bsmat{\left[ \begin{smallmatrix}}
\def\esmat{\end{smallmatrix} \right]}

\begin{document}

\begin{frontmatter}

\title{A unified method for optimal arbitrary pole placement}

\thanks[footnoteinfo]{Corresponding author: R.~Schmid.
 Earlier versions of this work appeared in \citep{SNNPa} and \cite{SNNPb}. This work was supported in part by the Australian Research Council under the grant FT120100604. }

\author[Melbourne]{Robert Schmid}\ead{rschmid@unimelb.edu.au}, 
\author[Perth]{Lorenzo Ntogramatzidis}\ead{L.Ntogramatzidis@curtin.edu.au}, 
\author[Exeter]{Thang Nguyen}\ead{T.Nguyen-Tien@exeter.ac.uk}, 
\author[SanDiego]{Amit Pandey}\ead{appandey@ucsd.edu} 

\address[Melbourne]{Department of Electrical and Electronic Engineering,
 University of Melbourne, Parkville, VIC 3010, Australia.}
\address[Perth]{Department of Mathematics and Statistics,
Curtin University, Perth, WA 6848, Australia.} %
\address[Exeter]{Department of Engineering, University of Exeter, UK.} %
\address[SanDiego]{Department of Mechanical and Aerospace Engineering, University of California, San Diego, USA. } %


\begin{abstract}

 We consider the classic problem of pole placement by state feedback. We offer an eigenstructure assignment algorithm to obtain a novel parametric form for the pole-placing feedback matrix that can deliver any set of desired closed-loop eigenvalues, with any desired multiplicities. This parametric formula is then exploited to introduce an unconstrained nonlinear optimisation algorithm to obtain a feedback matrix that delivers the desired pole placement with optimal robustness and minimum gain. Lastly we compare the performance of our method against several others from the recent literature.
\end{abstract}
\end{frontmatter}

\section{Introduction}

We consider the classic problem of repeated pole placement for linear time-invariant (LTI) systems in state space form
 \be
 \dot x(t) = A\,x(t)+B\,u(t),\;\quad \label{syseq1}
 \ee
where, for all $t \in \real $, $x(t) \in \real^n$ is the state and $u(t) \in \real^m$ is the control input.  We assume that $B$ has full column-rank, and that the pair $(A,B)$ is reachable. We let $\mathcal{L} = \{\lambda_1, \ldots, \lambda_\nu \}$ be a self-conjugate set of $\nu \le n$ complex numbers, with associated algebraic multiplicities $\mathcal{M} = \{m_1, \dots, m_\nu\}$ satisfying $m_1 + \dots + m_\nu= n$, and $m_i = m_j$ whenever $\lambda_i = \overline{\lambda}_j$. The problem of \textit{exact pole placement (EPP) by state feedback} is that of finding a real feedback matrix $F$ such that 
\be \label{eppeq}
(A+B\,F)\,X = X\,\Lambda,
\ee
where $\Lambda$ is a $n \times n$ Jordan matrix obtained from the eigenvalues of $\mathcal{L}$, including multiplicities given by $\mathcal{M}$, and $X$ is a matrix of closed-loop eigenvectors of unit length. The matrix $\Lambda$ can be expressed in the Jordan (complex) block  diagonal canonical form
\be \label{Jordan}
\Lambda = \blkdiag(\begin{array}{ccc} J(\lambda_1), & \cdots, & J(\lambda_\nu)
 \end{array}),
 \ee
 where each $J(\lambda_i)$ is a Jordan matrix for $\lambda_i$ of order $m_i$, and may be composed of up to $g_i$ mini-blocks
 \be \label{Jmini}
 J(\lambda_i) = \blkdiag(\begin{array}{ccc} J_1(\lambda_i), & \cdots, & J_{g_i}(\lambda_i) \end{array}),
 \ee
 where $1 \leq g_i \le m$. We use $\gP \defi \{p_{i,k}\,|\, 1 \leq i \leq \nu, 1 \leq k \leq g_i \}$ to denote the order of each Jordan mini-block $J_k(\lambda_i) $;  then
  $p_{i,k} = p_{j,k}$  whenever $\lambda_i = \overline{\lambda}_j$.
  When $(A,B)$ is reachable, arbitrary multiplicities of the closed-loop eigenvalues can be assigned by state feedback, but the possible mini-block orders of the Jordan structure of $A+BF$ are constrained by the \textit{controllability indices} \citep{Rosenbrock-70}.
 If $\gL$, $\gM$ and $\gP$ satisfy the conditions of the Rosenbrock theorem, we say that the triple $(\gL,\gM,\gP)$ defines an \textit{admissible Jordan structure} for $(A,B)$.

In order to consider optimal selections for the feedback matrix, it is important to have a parametric formula for the set of feedback  matrices that deliver the desired pole placement. In \citep{KNV} and \citep{SPN}  parametric forms are given for the case where
 $\Lambda$ is a diagonal matrix and the eigenstructure is non-defective; this requires $m_i \leq m$ for all $m_i \in \mathcal M$. Parameterisations that do not impose a constraint on the multiplicity of the eigenvalues to be assigned include \citep{BS} and \citep{FO83}; however these methods require the closed-loop eigenvalues to all be distinct from the open-loop ones.

The general case where $\mathcal L$ contains any desired closed-loop eigenvalues and multiplicities is considered in \citep{Chu} and \citep{RFBT}, where parametric formulae are provided for $F$ that use the eigenvector matrix $X$ as a parameter. Maximum generality in these parametric formulae has however been achieved at the expense of efficiency, as the square matrix $X$ has  $n^2$ free parameters. By contrast, methods \citep{KNV,FO83,BS,SPN} all employ parameter matrices with $mn$ free parameters.

The first aim of this paper is to offer a parameterisation for the pole-placing feedback matrix that combines the generality of \citep{Chu} and \citep{RFBT} with the efficiency of an $mn$-dimensional parameter matrix. We offer a parametric formula for all feedback  matrices $F$ solving (\ref{eppeq}) for any admissible $(\gL,\gM,\gP)$.
 For a given  parameter matrix $K$, we obtain the eigenvector matrix $X_{\scriptscriptstyle K}$ and feedback  matrix $F_{\scriptscriptstyle K}$ by building the Jordan chains from  eigenvectors  selected from the kernels of the matrix pencils $[A-\lambda_i\,I_n\;\; B]$, and thus avoid the need for matrix inversions, or the solution of Sylvester matrix equations. The parameterisation will be shown to be exhaustive of all feedback  matrices that assign the desired eigenstructure.

The second aim of the paper is to seek the solution to some optimal control problems. We firstly consider the \textit{robust exact pole placement problem} (REPP), which involves obtaining $F$ that renders the eigenvalues of $A + B\,F$ as insensitive to perturbations in $A$, $ B$ and $F $ as possible. Numerous results \citep{Ch} have appeared linking the sensitivity of the eigenvalues to various measures of the {\em condition number} of $X$. Another commonly used robustness measure is the {\em departure from normality} of the closed loop matrix $A + B\,F$.
For the case of diagonal $\Lambda$, there has been considerable literature on the REPP, including \citep{KNV,BN,TY,V2,RFBT,Chu,LCL,SPN}. Papers considering the REPP for the general case where $(\gL,\gM,\gP)$ defines an admissible Jordan structure include \citep{LTT} and \citep{RFBT}.

 A related optimal control problem is the \textit{minimum gain exact pole placement problem} (MGEPP), which involves solving the EPP problem and also obtaining the feedback matrix $F$ that has the least gain (smallest matrix  norm), which gives a measure of the control amplitude or energy required by the control action. Recent papers  addressing the MGEPP  with minimum Frobenius norm for $F$ include \citep{AE} and \citep{KU}.

 In this paper we utilise our parametric form for the matrices $X$ and $F$ that solve (\ref{eppeq}) to take a unified approach to the REPP and MGEPP problems, for any admissible Jordan structure. In our first method for the REPP, we seek the parameter matrix $K$ that minimises the Frobenius condition number of $X$. In our second approach to the REPP, we seek the parameter matrix that minimises the departure from normality of matrix $A+BF$. Next we address the MGEPP by seeking the parameter $K$ that minimises the Frobenius norm of $F$. Finally, we combine these approaches by introducing an objective function expressed as a weighted sum of robustness and gain measures, and use gradient iterative methods to seek a local minimum.

 The performance of the our algorithm will be compared against the methods of \citep{RFBT}, \citep{AE} and \citep{LCL} on a number of sample systems. We see that the methods introduced in this paper can achieve superior robustness while using less gain than all three of these alternative methods.

\section{Arbitrary pole placement}
\label{Moore}

Here we adapt the algorithm of \citep{KM} to obtain a simple parametric formula for the gain matrix $F$ that solves the exact pole placement problem for an admissible Jordan structure $(\gL,\gM,\gP)$, in terms of an arbitrary parameter matrix $K$ with  $mn$ free dimensions. We begin with some definitions.

 Given a self-conjugate set of $\nu$ complex numbers $ \{\lambda_1,\ldots,\lambda_{\nu}\}$ containing  $\sigma$ complex conjugate pairs, we say that the set is  {\em $\sigma$-conformably ordered} if the first $2\,\sigma$ values are complex while the remaining are real, and for all odd $i \le 2\,\sigma$ we have $\lambda_{i+1}=\overline{\lambda}_i$. For example, the set $\{10\,j,-10\,j,2+2\,j,2-2\,j,7\}$
 is $2$-conformably ordered. For simplicity we shall assume in the following that $\gL$ is $\sigma$-conformably ordered.

     If  $M$ is a complex  matrix partitioned into $\nu$ column matrices $M= [M_1  \ldots   M_\nu]$, we say that $M$ is {\em $\sigma$-conformably ordered} if the first $2\,\sigma$ column matrices  of $M$ are complex while the remaining are real, and for all odd $i \le 2\,\sigma$ we have $M_{i+1}=\overline{M}_i$. For a $\sigma$-conformably ordered complex  matrix $M$, we  define a  real matrix $Re(M)$ composed of $\nu$ column matrices of the same dimensions as  those of $M$ thus: for each  odd $i \in \{1,\dots,2\sigma\}$, the $i$-th and $i+1$-st column matrices   of  $Re(M)$ are $\frac{1}{2}(M_i+M_{i+1})$ and $\frac{1}{2j}(M_i-M_{i+1})$ respectively,  while for  $i \in \{2\sigma +1, \dots,\nu\}$, the column matrices of $Re(M)$ are the same as the corresponding column matrices  of $M$. For any real or complex  matrix $X$ with  $n+m$ rows, we define matrices $\overp(X)$ and  $\undp(X)$  by taking the first $n$ and last $m$ rows of $X$, respectively.
For each $i \in \{1,\dots,\nu\}$, we define the matrix pencil
 \be \label{Slambda}
 S(\lambda_i)\defi \bmat{cccc} A-\lambda_i \,I_n &&& B\emat.
 \ee
   We use  $N_i$ to  denote an orthonormal basis matrix for the kernel of $S(\lambda_i)$. If $\lambda_{i+1}=\overline{\lambda}_i$, then $N_{i+1}=\overline{N}_i$.  Since each $S(\lambda_i)$ is $n \times (n+m)$ and $(A,B)$ is reachable, each kernel  has  dimension $m$.
 We let
 \be \label{Mlambda}
 M_i\defi \bmat{cccc} A-\lambda_i\, I_n &&& B\emat^\dagger,
 \ee
 where ${}^\dagger$ indicates the Moore-Penrose pseudo-inverse. For any matrix $X$ we use $X(l)$ to denote the $l$-th column of $X$.

We say that a matrix $K$ is a
 {\em compatible parameter matrix} for $(\gL,\gM,\gP)$, if $ K \defi \blkdiag\{ K_1, \dots, K_\nu\}$, where each $K_i$ has dimension $m \times m_i$, and for each $i \ge 2\,\sigma$, $K_i$ is a real matrix, and for all odd $i \leq 2\,\sigma$, we have $K_{i+1} =\overline{K}_{i}$. Then each $K_i$ matrix may be partitioned as
 \be
 \label{Kieq}
 K_i = \bmat{cccc} K_{i,1} \ & K_{i,2} \ & \dots \ & K_{i,g_i}\emat,
 \ee
where each $K_{i,k}$ has dimension $m \times p_{i,k}$. For  $i \in \{1,\ldots,\nu\}$ and  $k \in \{1,\dots, g_i\}$ we build vector chains of length $p_{i,k}$ as
\bea
h_{i,k}(1) & \;=\; & N_i \,K_{i,k}(1), \label{eq32} \\
h_{i,k}(2) & \;=\; & M_i \,\overp \{ h_{i,k}(1)\} + N_i \,K_{i,k}(2), \label{eq33} \\
 & \vdots & \nn \\
h_{i,k}(p_{i,k}) & \;=\; & M_i \, \overp \{h_{i,k}(p_{i,k}-1) \}+ N_i\, K_{i,k}(p_{i,k}). \label{eq34}
\eea
From these column vectors we construct the matrices
\be
H_{i,k} \defi [ h_{i,k}(1)\  \dots \ h_{i,k}(p_{i,k})] \label{Hik}
\ee
of dimension $ (n + m) \times p_{i,k}$,  and
\be \label{Hidef}
H_i \defi [H_{i,1} \ \dots \ H_{i,g_i} ], \ \
H_K \defi  [ H_1 \ \dots \ H_\nu ], \ \ X_{\scriptscriptstyle K}  \defi  \overp \{H_{\scriptscriptstyle K}\}
\ee
of dimension  $ (n+m) \times m_i$, $ (n+m) \times n$ and $n \times n$, respectively. Note that $H_K$ is $\sigma$-conformably ordered,  and  hence we  may define  real matrices
%
\be
 \quad
V_{\scriptscriptstyle K}  \defi \overp \{Re(H_{\scriptscriptstyle K})\}, \quad W_{\scriptscriptstyle K} \defi \undp \{Re(H_{\scriptscriptstyle K})\}  \label{Veq} \ee
of dimensions $n \times n$ and $m \times n$, respectively.
 We are now ready to present the main result of this paper.
 \begin{theorem} \label{Mooreprop}
 For almost all choices of the compatible  parameter matrix $K$, the matrix $V_{\scriptscriptstyle K}$ in (\ref{Veq}) is invertible. The set of all real feedback matrices $F$ such that $A+B\,F$ has Jordan structure given by $(\gL, \gM,\gP)$ is parameterised in $K$ as
 \be \label{Feq}
 F_{\scriptscriptstyle K} = W_{\scriptscriptstyle K}\,V_{\scriptscriptstyle K}^{-1}.
 \ee

\end{theorem}

\proof 
 Firstly we let $K$ be any  compatible parameter matrix yielding invertible $V_{\scriptscriptstyle K}$ and $W_{\scriptscriptstyle K}$ in  (\ref{Veq})  and  $F_{\scriptscriptstyle K}$ in (\ref{Feq}). We  prove that the closed-loop matrix $A+BF_{\scriptscriptstyle K}$ has the required eigenstructure.
  $V_K$ and $W_K$ may be partitioned as
\be
V_{\scriptscriptstyle K}  =  [V_1\;\;\ldots \;\; V_{\nu}], \quad  W_{\scriptscriptstyle K}  =  [W_1\;\; \ldots\;\; W_{\nu}], \label{Vkeq}
\ee
where, for each $i \in \{1,\ldots,\nu\}$,  $V_i$ and $W_i$ have $m_i$ columns. Let $H_{i,k}$  in (\ref{Hik}) be partitioned  as
 \be
 H_{i,k} = \left[ \begin{array}{ccc}
 v_{i,k}^\prime(1) \ & \dots \ & v_{i,k}^\prime(p_{i,k}) \\
 w_{i,k}^\prime(1) \  & \dots \ & w_{i,k}^\prime(p_{i,k}) \end{array} \right],
 \ee
 where, for each $k \in \{1,\dots,g_i\}$,  the column vectors satisfy by construction
 \bq
(A -\lambda_i\,I_n) v_{i,k}^\prime(1) + B w_{i,k}^\prime(1) & = & 0, \label{ML1}\\
(A -\lambda_i\,I_n) v_{i,k}^\prime(2) + B w_{i,k}^\prime(2) & = & v_{i,k}^\prime(1), \label{ML2} \\
 & \vdots & \nonumber \\
(A -\lambda_i\,I_n) v_{i,k}^\prime(p_{i,k}) + B w_{i,k}^\prime(p_{i,k}) & = & v_{i,k}^\prime(p_{i,k}-1). \quad \label{ML3} \eq
Define for each $i \in \{1,\ldots,\nu\}$ and $k \in \{1,\dots,g_i\}$,
\be
V_{i,k}^\prime  =  [ v_{i,k}^\prime(1) \   \dots \ v_{i,k}^\prime(p_{i,k})], \quad
W_{i,k}^\prime  =  [ w_{i,k}^\prime(1) \  \dots \  w_{i,k}^\prime(p_{i,k})],
\ee
and next define, for each $i \in \{1,\ldots,\nu\}$,  $V_{i}^\prime = [ V_{i,1}^\prime \ \dots \  V_{i,g_i}^\prime ]$ and $W_{i}^\prime = [ W_{i,1}^\prime \ \dots \ W_{i,g_i}^\prime ]$.  As  $K$ is a compatible parameter matrix,  we  have, for all odd $i \in \{1,\ldots,2\,\sigma\}$, $V_{i+1}^\prime=\overline{V}_{i}^\prime$ and $W_{i+1}^\prime=\overline{W}_{i}^\prime$.  Finally, introduce
$U_i \defi \frac{1}{2}\bsmat I_{m_i} & -j\, I_{m_i} \\[1mm] I_{m_i} & j\,I_{m_i}\esmat$. Then for each odd $i \in \{1,\ldots,2\,\sigma\}$, we have $[ V_i^\prime \;\; V_{i+1}^\prime ] \,U_i= [ V_i \;\; V_{i+1} ]$ and
$[ W_i^\prime \;\; W_{i+1}^\prime ] \,U_i= [ W_i \;\; W_{i+1} ]$,  and for each $i \in \{2\sigma+1, \dots, \nu\}$, we  have
$ V_i^\prime =  V_i $ and $ W_i^\prime =  W_i $.
Since $F_{\scriptscriptstyle K}\,V_{\scriptscriptstyle K}=W_{\scriptscriptstyle K}$, then
$F_{\scriptscriptstyle K}\,[V_i^\prime \;\; V_{i+1}^\prime]=[W_i^\prime \;\; W_{i+1}^\prime]$
for all odd $i \in \{1,\ldots,2\,\sigma\}$ and $F_{\scriptscriptstyle K}\,V_i=W_i$ for all $i \in \{2\,\sigma+1,\ldots,\nu\}$.
 Hence, for each odd $i\in\{1,\ldots,2\,\sigma\}$,  we  have
  \be
 (A+B\,F_{\scriptscriptstyle K})[\, V_i^\prime \;\; V_{i+1}^\prime \,]
 =[\, V_i^\prime \;\; V_{i+1}^\prime \,]\diag\{J(\lambda_i), J(\lambda_{i+1})\},
 \ee
 and   for all $i \in \{2\,\sigma+1,\ldots,\nu\}$, we have
$(A+B\,F_{\scriptscriptstyle K})V_i = V_i\, J(\lambda_i)$.  Thus
$(A+B\,F_{\scriptscriptstyle K})\,X_{\scriptscriptstyle K} = X_{\scriptscriptstyle K}\,\Lambda$,
where $X_{\scriptscriptstyle K} = [ \,V_1^\prime \;\; \ldots V_{\nu}^\prime\,]$ and $\Lambda$ is as in (\ref{Jordan}), as required.

 In order to prove that the parameterisation is exhaustive, we consider a feedback matrix $F$ such that the eigenstructure of $A+B\,F$ is given  by ($\gL,\gM, \gP)$, and  show  there exists a compatible parameter matrix  $K$ such that  matrices $V_{\scriptscriptstyle K}$ and $W_{\scriptscriptstyle K}$ can be constructed in (\ref{Veq}), with $V_{\scriptscriptstyle K}$ invertible and $F=W_{\scriptscriptstyle K}\,V_{\scriptscriptstyle K}^{-1}$.
   From   (\ref{Jordan})-(\ref{Jmini}),  $\Lambda$ can be written  as
\be
\Lambda = \blkdiag
(J_{1}(\lambda_1), \ldots,
J_{g_1}(\lambda_1),\ldots,
J_{1}(\lambda_\nu),\ldots,J_{g_\nu}(\lambda_\nu)).  \nn
\ee
Hence there  exists  an invertible matrix $T$ satisfying
$(A+B\,F) \,T= T \Lambda$. Let us partition $X$ and $Y$ conformably with the corresponding Jordan mini-blocks that they multiply, i.e.,
\be
 \bmat{cc} A  &  B  \emat \bmat{ccc}  X_{1,1} &  \ldots &   X_{\nu,g_\nu} \\  Y_{1,1} &  \ldots  & Y_{\nu,g_\nu} \emat = \left[ \begin{array}{ccc}  X_{1,1}\,J_1(\lambda_1) \ldots    &  X_{\nu,g_\nu}\,J_{g_\nu}(\lambda_\nu) \end{array} \right].  \nn
\ee
For  $i \in \{1,\ldots,\nu\}$  and  $k \in \{1,\dots,g_i\}$, the generic term is
\be
\label{alpha}
\bmat{cc} A & B \emat \bmat{c} X_{i,k} \\ Y_{i,k} \emat = X_{i,k} \, J_k(\lambda_i).
\ee
 First consider the case in which $\lambda_i$ is real. Partitioning $X_{i,k}=[v_{i,k}(1)\;\; \ldots\; v_{i,k}(p_{i,k})]$ and $Y_{i,k}=[w_{i,k}(1)\;\; \ldots\; w_{i,k}(p_{i,k})]$, we can write (\ref{alpha}) as
\be
 \bmat{cc} A  & B \emat \ns \bmat{ccc} v_{i,k}(1) & \ns \dots \ns &  v_{i,k}(p_{i,k}) \\ w_{i,k}(1) &\ns \dots  \ns &  w_{i,k}(p_{i,k}) \emat
  =\bmat{ccc} v_{i,k}(1) & \ns \dots \ns & v_{i,k}(p_{i,k})\emat  J_k(\lambda_i),  \nn
 \ee
 which yields
\bq
A\,v_{i,k}(1)+B\,w_{i,k}(1) &=&  v_{i,k}(1) \,\lambda_i \label{from1} \\
A\,v_{i,k}(2)+B\,w_{i,k}(2) &=& v_{i,k}(1)+\lambda_i\,v_{i,k}(2) \label{from2} \\
 \ns&\ns \vdots \ns&\ns  \nn \\
A\,v_{i,k}(p_{i,k}) +B\,w_{i,k}(p_{i,k})  &=& v_{i,k}(p_{i,k}-1)+\lambda_i\,v_{i,k}(p_{i,k}) \label{fromk}
\eq
We denote $h_{i,k}(l)= \bsmat v_{i,k}(l) \\[1mm] w_{i,k}(l)\esmat$.
From (\ref{from1}) we see that  $h_{i,k}(1) \in \ker (S(\lambda_i))$  and  hence there exists $K_{i,k}(1)$ satisfying (\ref{eq32}).  Moreover, from (\ref{from2}) we find
$[\begin{array}{cc}  A-\lambda_i\,I_n & B \end{array} ]h_{i,k}(2) =v_{i,k}(1)$,
which implies that there exists $K_{i,j}(2)$ satisfying (\ref{eq33}).
Repeating this procedure for all $l \in \{1,\ldots,p_{i,k}\}$, we find the parameters $K_{i,k}(1), \ldots, K_{i,k}({p_{i,k}})$ which satisfy (\ref{eq32})-(\ref{eq34}).  This procedure can be carried out for all real Jordan mini-blocks. Consider now the case of a real mini-block associated with a complex conjugate eigenvalue $\lambda_i=\sigma_i+j\,\omega_i$.
 For brevity we shall assume $p_{i,k}=2$. 
 Thus, (\ref{from1}) becomes
\beann
&&\bmat{cc}  A  &  B  \emat  \bmat{cccc} v_{i,k}(1) & v_{i,k}(2) & v_{i+1,k}(1) & v_{i+1,k}(2) \\  w_{i,k}(1) & w_{i,k}(2) & w_{i+1,k}(1) & w_{i+1,k}(2) \emat \\
&&
=\bmat{cccc} v_{i,k}(1) & v_{i,k}(2) & v_{i+1,k}(1) & v_{i+1,k}(2) \emat\!\! \bmat{cc|cc}
\sigma_i  &  \omega_i  &  1  &  0  \\
 -\omega_i  &  \sigma_i  &  0  &  1  \\
 \hline
0  &  0  &  \sigma_i  &  \omega_i  \\
 0  &  0  &  -\omega_i  &  \sigma_i \emat\!\!,
\eeann
which can be re-written as
\beann
&&\bmat{cc}  A  &  B  \emat \bmat{cc}  v_{i,k}(1)  +  j\,v_{i,k}(2)  &  v_{i+1,k}(1)  +  j\,v_{i,k}(2)  \\  w_{i,k}(1)  +  j\,w_{i,k}(2)  &  w_{i+1,k}(1)  +  j\,w_{i,k}(2)\emat  \\
&&=\!  \bmat{cc} \!\! v_{i,k}(1) \! +\!  j\,v_{i,k}(2)  &  v_{i+1,k}(1)\!  + \! j\,v_{i,k}(2)   \\ \!\! 0  &  0
\emat\!\!\!  \bmat{cc} \!\! \sigma_i \! + \! j\,\omega_i \! \!&\!  1  \!\!\! \\ \!\! 0  &\!\! \sigma_i  \!+\!  j\,\omega_i \!\! \emat\!\!,
\eeann
and the arguments above can be utilised after a re-labeling of the vectors.

Lastly we show that  $V_{\scriptscriptstyle K}$ is invertible for almost all choices of the parameter matrix $K$.
For each $i \in \{1,\dots,\nu\}$,  we  may express  the  orthonormal  basis $N_i$ for $\ker (S(\lambda_i))$ as
 $  N_i = [ h_{i,1} \ \dots \  h_{i,m} ]$. For each $k \in \{1, \dots, g_i\}$ we construct
\bq
h_{i,k}(1) & = & h_{i,k}  \\
h_{i,k}(2) & = & M_i\, h_{i,k}(1)  \\
 & \vdots & \nn \\
h_{i,k}(p_{i,k}) & = & M_i\, h_{i,k}(p_{i,k}-1)
\eq
and combining these we obtain
\be
H_{i,k}  =   [\, h_{i,k}(1)\,  \, \dots\, \, h_{i,k}(p_{i,k})].
\ee
Lastly we  obtain matrices $H_i $ and $H$ as  in (\ref{Hidef}), and $V$ as in (\ref{Veq}).
Then we must have $\rank(V)=n$, else no parameter matrix $K$ would exist to yield a real  feedback  matrix $F_{\scriptscriptstyle K}$ in (\ref{Feq}) that delivers the desired closed-loop eigenstructure. This contradicts the assumption that $(A,B)$ is reachable.

Next let $K$ be any compatible parameter matrix for $(\gL,\gM,\gP)$, let $V_{\scriptscriptstyle K} = \overp\{Re(H_{\scriptscriptstyle K})\}$ and assume
 $V_{\scriptscriptstyle K}$ is singular. Then $X_{\scriptscriptstyle K}$ in (\ref{Hidef}) is also  singular, i.e. $\rank(X_{\scriptscriptstyle K}) \leq n-1$.
 Without loss of generality, assume the first column of $X_K$ is linearly dependent upon the remaining ones. Then there exist a $\sigma$-conformalby ordered set of $n$ coefficient vectors
 $\alpha_{i,k,l}$, not all equal to zero, for which
 \bq
\overp\{h_{1,1}(1)K_{1,1}(1)\}
& = &  \sum_{l=2}^{p_{1,1}}  \alpha_{1,1,l}\ \overp\{h_{1,1}(l)\} \nn \\
&  & + \sum_{k=2}^{g_1} \sum_{l=1}^{p_{1,k}} \ \alpha_{1,k,l} \overp\{h_{1,k}(l)\} \nn \\
&  & + \sum_{i =2}^{\nu} \sum_{k=1}^{g_i} \sum_{l=1}^{p_{i,k}}\  \alpha_{i,k,l} \overp\{h_{i,k}(l)\}  \nn
\eq
 This implies that $\rank(X_{\scriptscriptstyle K})=n$ may fail only when $K_{1,1}(1) $ lies on an $(m-1)$-dimensional hyperplane in the $m$-dimensional parameter space. Thus the set of compatible parameter matrices $K$ that can lead to a loss of rank in $X_{\scriptscriptstyle K}$,  and  hence  $V_{\scriptscriptstyle K}$,  is given by the union of at most  $n$ hyperplanes of dimension at most $nm-1$ in  the $nm$-dimensional parameter space. Since hyperplanes have zero Lebesgue measure, the set of parameter matrices $K$ leading to singular $V_{\scriptscriptstyle K}$ has zero Lebesgue measure.
\endproof

 The above formulation takes its inspiration from the proof of Proposition 1 in \citep{KM}, and hence we shall refer to (\ref{Feq}) as the {\it Klein-Moore parametric form} for $F$.

\section{Optimal pole placement methods} \label{secopt}

We firstly present some classic results on eigenvalue sensitivity.  Let $A$ and $X$ be such that $A = XJX^{-1}$, where $J$ is the Jordan form of $A$, and let $A^\prime = A + H$. Then, for each eigenvalue $\lambda^\prime $ of $A^\prime$, there exists an eigenvalue $\lambda$ of $A$ such that
 \be
 \frac{|\lambda - \lambda^\prime |}{(1 + |\lambda - \lambda^\prime |)^{l-1}} \leq \kappa_2(X)\|H\|_2,
 \ee
 where $l$ is the size of the largest Jordan mini-block associated with $\lambda$,  and $\kappa_2(X)\defi \|X\|_2\|X^{-1}\|_2$ is the spectral  condition number of $X$  \citep{Ch}.
 As the Frobenius condition number $\kappa_{\textsc{fro}}(X)= \|X\|_{\textsc{fro}}\|X^{-1}\|_{\textsc{fro}}$ satisfies $\kappa_2(X) \leq \kappa_{\textsc{fro}}(X)$ and is differentiable, it is often used as a robustness measure in conjunction with gradient search methods.

A second widely used robustness measure is the {\em departure from normality} of the matrix $A$, which is defined as follows \citep{SS}:
 Let $U$ be any unitary matrix such that $U^\tra\,A\,U$ is upper triangular, then $U^\tra\,A\,U = D + R$, for some diagonal matrix $D$ and strictly upper triangular matrix $R$. The Frobenius departure from normality of $A$ is then $\delta_{\textsc{fro}}(A)\defi \|R\|_{\textsc{fro}}$.

 Our Method 1 simultaneously addresses the REPP and MGEPP by using the weighted objective function
 \be
 f(K)=\alpha \kappa_{\textsc{fro}}(V_{\scriptscriptstyle K}) +(1 - \alpha)\|F_{\scriptscriptstyle K}\|_{\textsc{fro}},
\ee
where $K$ is a compatible parameter  matrix and $V_{\scriptscriptstyle K}$ and $F_{\scriptscriptstyle K}$ are obtained from (\ref{Veq}) and (\ref{Feq}). Finding $K$ to minimise $f$ presents an unconstrained nonconvex optimisation problem. For efficient computation \citep{BN}, showed we can use the equivalent objective function
\be
 f_1(K)=\alpha(\|V_{\scriptscriptstyle K}\|_{\textsc{fro}}^2+ \|V_{\scriptscriptstyle K}^{-1}\|_{\textsc{fro}}^2) +(1 - \alpha)\|F_{\scriptscriptstyle K}\|_{\textsc{fro}}^2, \label{mgrepp1}
\ee
Here, $\alpha$ is a weighting factor, with $0 \leq \alpha \leq 1$. The limiting cases $\alpha = 0$ and $\alpha =1$ define the MGEPP and REPP problems, respectively.

 Our Method 2 uses the weighted objective function
 \be
 f_2(K)=\alpha \delta_{\textsc{fro}}^2(A+B\,F_{\scriptscriptstyle K}) +(1-\alpha)\|F_{\scriptscriptstyle K} \|_{\textsc{fro}}^2. \label{mgrepp2}
\ee
Finding $K$ to minimise $f_2$ again presents an unconstrained nonconvex optimisation problem.
Expressions for the derivatives of  $H_{\scriptscriptstyle K}$, $\|V_{\scriptscriptstyle K}\|_{\textsc{fro}}$,  and $\|V_{\scriptscriptstyle K}^{-1} \|_{\textsc{fro}}$   were given in \cite{SNNPb}; from these, gradient search methods can be used to seek local minima for $f_1$ and $f_2$. The results are contingent upon the initial choice of the parameter matrix $K$.
%

\section{Performance comparisons}
In this section, we compare the performance of our algorithm with the methods given in the recent papers by \citep{RFBT}, \citep{AE} and \citep{LCL}. In \citep{BN} a collection of benchmark systems were introduced that have been investigated over the years by many authors. To compare our performance against the method of \citep{RFBT}, we used the matrices $(A,B)$ from these examples, but in order to compare their performance for defective pole assignment, we assigned all the closed-loop eigenvalues to zero. In each case we assigned Jordan blocks of sizes equal to the controllability indices. Using the toolbox \textit{rfbt} to implement the method of \citep{RFBT} that we created for our earlier computational survey in \citep{SPN}, we obtained the matrices $F$ and $X$ delivered by this method, for each of the 11 sample systems. We also implemented our own method on these systems. The results are shown in Table \ref{Surv1}.

\begin{table}[tbh!]
\begin{center}
\caption{Byers and Nash examples} \label{Surv1}
\begin{tabular}{||c||c|c||c|c||}
 \hline
 Example & \multicolumn{2}{c||}{Ait Rami \em{et al} } & \multicolumn{2}{c||}{Our Method}
 \\ \cline{2-5}
 & $\kappa_{\textsc{fro}}(X)$ & $\|F\|_{\textsc{fro}}$ & $\kappa_{\textsc{fro}}(X)$ & $\|F\|_{\textsc{fro}}$ \\
 \hline
 1 & 16.73 & 3.102 & 16.73 & 3.102 \\[-1mm]
 2 & 54.43 & 645.5 & 51.11 & 289.5 \\[-1mm]

 3 & 7.188 & 2.225 & 7.188 & 2.225 \\[-1mm]

 4 & 11.49 & 7.145 & 11.49 & 7.043 \\[-1mm]

 5 & 29.99 & 186.8 & 28.39 & 138.0 \\[-1mm]

 6 & 113.4 & 8.167 & 113.4 & 7.880 \\[-1mm]

 7 & 16.84 & 595.9 & 17.33 & 596.1 \\[-1mm]

 8 & 4.000 & 10.07 & 4.000 & 9.230 \\[-1mm]

 9 & 85.68 & 22,610 & 85.65 & 22,610 \\[-1mm]

 10 & 30.33 & 29.74 & 30.33 & 29.74 \\[-1mm]

 11 & 4,579 & 5,025 & 4,501 & 5,025 \\[0mm]
\hline
 \end{tabular}
\end{center}
\end{table}

 Comparing the robust conditioning performance of the two methods, we see little difference between the methods. However, when we compare the matrix gain used to achieve this eigenstructure we observe that our method was able to use less gain in 5 of the sample systems, and in two cases (System 2 and 5) the reduction in gain was very considerable. The results are in agreement with the findings of the survey in \citep{SPN}, which considered sample systems with non-defective eigenstructure and found that our method could achieve comparable robust conditioning with that of \cite{RFBT}, but with reduced gain.
 %

To compare our performance against that of \citep{AE}, we considered the 5 example systems introduced in that paper. Among these, the first example system assigned all the poles to zero, and hence requires a defective closed-loop eigenstructure. The other four sample systems all involve distinct eigenvalues. The results are shown in Table \ref{Surv2}. The results have been constructed using the feedback  matrices provided by \citep{AE}.

\begin{table}[tbh!]
\caption{Ataei and Enshaee examples} \label{Surv2}
\begin{center}
\begin{tabular}{||c||c|c||c|c||}
 \hline
 Example & \multicolumn{2}{c||}{Ataei and Enshaee} & \multicolumn{2}{c||}{Our Method}
 \\ \cline{2-5}
 & $\kappa_{\textsc{fro}}(X)$ & $\|F\|_{\textsc{fro}}$ & $\kappa_{\textsc{fro}}(X)$ & $\|F\|_{\textsc{fro}}$ \\
 \hline
 1 & 321.4 & 1.295 & 4.444 & 1.295 \\[-1mm]

 2 & 290.5 & 3.970 & 278.6 & 3.844 \\[-1mm]

 3 & 7.895 & 1.311 & 6.515 & 1.304 \\[-1mm]

 4 & 3.873 & 4.243 & 4.353 & 4.072 \\[-1mm]

 5 & 26.01 & 4.748 & 21.56 & 4.662 \\
\hline
 \end{tabular}
\end{center}
\end{table}

 The results show that our method achieved the desired eigenstructure with equal or slightly less gain than that of \citep{AE}. In all but one of the samples, our method also achieved a more robust eigenstructure, especially in Example 1, which has the defective eigenstructure.

 Lastly, we consider Example 1 in \citep{LCL}. The four desired closed loop poles are all distinct in this example. The method of \citep{LCL} considers the problem of minimising the Frobenius norm of the feedback  matrix and the minimisation of the departure from normality measure. The authors obtained a feedback  $F$ yielding $\delta_{\textsc{fro}}(A +BF) = 20.67$, and an alternative  matrix $F$ that delivers the desired pole placement with gain $\|F\|_{\textsc{fro}} = 6.049$.

 Applying Method 2 with $\alpha = 1$ we obtained a feedback  matrix $F$ yielding $\delta_{\textsc{fro}}(A +BF) = 18.52$, and by using $\alpha = 0$, we obtained $F$ such that $\|F\|_{\textsc{fro}} = 3.826$, indicating that our method can achieve the desired pole placement with either smaller  departure from normality measure, or less gain, than the method of \citep{LCL}, as required.

\section{Conclusion}
We have introduced a novel parametric form for the feedback  matrix that solves the classic problem of exact pole placement with any desired eigenstructure. The parametric form was used to take a unified approach to a variety of optimal pole placement problems. The effectiveness of the method has been compared against several recent alternative methods from the literature, and was shown in several examples to achieve the desired pole placement with either superior robustness or smaller gain than the other methods surveyed.


\bibliographystyle{plainnat}

\begin{thebibliography}{100}


\bibitem[Ait Rami {\em et al}(2009)]{RFBT}
M.~Ait Rami, S.E.~Faiz, A.~Benzaouia, and
F.~Tadeo, Robust Exact Pole Placement via an LMI-Based Algorithm, {\it IEEE Transactions on Automatic Control}, 54, 394--398, 2009.


\bibitem[Ataei and Enshaee(2011)]{AE}
M.~Ataei and A.~Enshaee, Eigenvalue assignment by minimal state-feedback gain in LTI multivariable systems, {\it International Journal of Control}, 84, 1956-–1964, 2011.


\bibitem[Bhattacharyya and de Souza(1982)]{BS}
S.P.~Bhattacharyya and E.~de Souza, Pole assignment via Sylvester equation, {\it Systems \& Control Letters}, 1, 261--263, 1982.


\bibitem[Byers and Nash(1989)]{BN} R.~Byers and S.G.~Nash, Approaches to robust pole assignment, {\it International Journal of Control}, 49, 97-–117, 1989.

\bibitem[Chatelin(1993)]{Ch} F.~Chatelin, Eigenvalues of Matrices, John Wiley and Sons, 1993.

\bibitem[Chu(2007)]{Chu}
E.~Chu, Pole assignment via the Schur form, {\it Systems \& Control Letters} 56, 303-–314, 2007.


\bibitem[Fahmy and O'Reilly(1983)]{FO83}
M.M.~Fahmy and J.~O'Reilly, Eigenstructure Assignment in Linear Multivariable Systems - A Parametric Solution, {\it IEEE Transactions on Automatic Control}, 28,
 990--994, 1983.
%
%
%
%


\bibitem[Kautsky {\em et al}(1985)]{KNV}
J.~Kautsky, N.K.~Nichols and P.~van Dooren, Robust Pole Assignment in Linear State Feedback, {\it International Journal of Control},
 41, 1129--1155, 1985.

\bibitem[Klein and Moore(1977)]{KM}
G.~Klein and B.C.~Moore, Eigenvalue-Generalized Eigenvector Assignment with State Feedback, {\it IEEE Transactions on Automatic Control}, 22,
 141--142, 1977.

\bibitem[Lam {\em et al}(1997)]{LTT} J.~Lam, H.K.~Tam and N.K.~Tsing, Robust deadbeat regulation, {\it International Journal of Control},
 67, 587--602, 1997.

\bibitem[Li {\em et al}(2011)]{LCL}
T.~Li, E.~Chu and W.W.~Lin, Robust Pole Assignment for Ordinary
and Descriptor Systems via the Schur Form, in {\it Numerical Linear Algebra in Signals, Systems and Control}, edited by P. van Dooren,
Lecture Notes in Electrical Engineering (80), Springer 2011.

\bibitem[Kochetkov and Utkin(2014)]{KU}
S.A.~Kochetkov and V.A.~Utkin, Minimizing the Feedback Matrix Norm in Modal Control Problems, {\it Automation and Remote Control}, 75, 234–-262, 2014.
%

\bibitem[Rosenbrock(1970)]{Rosenbrock-70}
H.H.~Rosenbrock, State-Space and Multivariable Theory. New York: Wiley, 1970.

\bibitem[Schmid {\em et al}(2013a)]{SNNPa}
 R.~Schmid, L.~Ntogramatzidis, T.~Nguyen and A.~Pandey, Arbitrary pole placement with minimum gain, {\it Proceedings 21st Mediterranean Conference on Control and Automation}, Crete, 2013a.


\bibitem[Schmid {\em et al}(2013b)]{SNNPb}
R.~Schmid, L.~Ntogramatzidis, T.~Nguyen and A.~Pandey, Robust repeated pole placement, {\it Proceedings 3rd
 IEEE Australian Control Conference}, Melbourne, 2013b.


\bibitem[Schmid {\em et al}(2014)]{SPN}
R.~Schmid, A.~Pandey and T.~Nguyen, Robust Pole Placement With Moore's Algorithm,
 {\it IEEE Transactions on Automatic Control}, 59, 500--505, 2014.

\bibitem[Stewart and Sun(1990)]{SS}
G.W.~Stewart and J.G.~Sun, Matrix Perturbation Theory, Academic press, 1990.

\bibitem[Tam and Lam(1997)]{TL}
H.K.~Tam and J.~Lam,
Newton's approach to gain-controlled robust pole placement, {\it IEE Proc.-Control Theory Applications}, 144, 439--446, 1997.

\bibitem[Tits and Yang(1996)]{TY}
A L.~Tits and Y.~Yang, Globally Convergent Algorithms for Robust Pole
Assignment by State Feedback, {\it IEEE Transactions on Automatic Control}, 41, 1432-–1452, 1996.

\bibitem[Varga(2000)]{V2}
A.~Varga, Robust Pole Assignment via Sylvester Equation Based State Feedback Parametrization,
 {\it Proceedings IEEE International Symposium on Computer-Aided Control System Design}, Anchorage, USA, 2000.

\end{thebibliography}

\end{document}